\documentclass[11pt]{article}

\usepackage{amsmath,amssymb,amsthm,verbatim,amscd,eucal}
\usepackage{pictexwd}
\usepackage{latexsym}
\usepackage{graphicx}
\usepackage[usenames]{color}
\usepackage{pstricks,pst-plot,pst-node,pst-text,pst-3d}

\numberwithin{equation}{section}
\newtheorem{thm}[equation]{Theorem}
\newtheorem{prop}[equation]{Proposition}

\newtheoremstyle{efronremark}
{6pt}{6pt}{}{}{\itshape}{\quad}{ }{\thmnote{#3}}

\theoremstyle{efronremark}   \newtheorem*{eremark}{}

\newcommand{\ZZ}{\mathbb{Z}}
\newcommand{\CC}{\mathbb{C}}
\newcommand{\C}{\mathsf{C}}
\newcommand{\N}{\mathsf{N}}

\newcommand{\LL}{\mathsf{L}}

\newcommand{\M}{\mathsf{M}}
\newcommand{\PS}{\mathsf{P}}
\newcommand{\qq}{\mathsf{q}}

\newcommand{\RR}{\mathsf{R}}

\newcommand{\UU}{\mathsf{U}}
\newcommand{\VV}{\mathsf{V}}

\newcommand{\g}{\mathfrak{g}}
\newcommand{\gl}{\mathfrak{gl}(1|1)}
\newcommand{\uqgl}{\UU_{\mathsf q}(\gl)}

\begin{document}

\title{Planar Rook Algebras and Tensor Representations of $\gl$}

\author{\textsc{Georgia Benkart} \\
\textit{Department of Mathematics} \\
\textit{University of Wisconsin - Madison}  \\
\textit{Madison, WI 53706, USA}\\
\texttt{benkart@math.wisc.edu}
\and
\textsc{Dongho Moon} \\
\textit{Department of Mathematics}\\
\textit{Sejong University}\\
\textit{Seoul 133-747 Korea (ROK)}\\
\texttt{dhmoon@sejong.ac.kr}}
\date{January 11,  2012}
\maketitle

\begin{abstract}    We establish a connection between planar rook algebras
and tensor representations $\VV^{\otimes k}$ of the natural two-dimensional representation $\VV$  of the general linear Lie superalgebra $\gl$.  In particular, we show that the centralizer algebra $\mathsf{End}_{\gl}(\VV^{\otimes k})$ is the planar rook algebra $\CC \mathsf{P}_{k-1}$ for all $k \geq 1$,  and we exhibit an explicit decomposition of
 $\VV^{\otimes k}$  into irreducible $\gl$-modules.    We obtain similar results for the quantum
 enveloping algebra $\UU_\qq(\gl)$ and its natural two-dimensional module $\VV_\qq$.
 \bigskip

\begin{eremark}[Keywords:] planar rook algebra, general linear Lie superalgebra, $\gl$-representation, quantum enveloping algebra $\UU_\qq(\gl)$
\end{eremark}
\begin{eremark}[2010 Mathematics Subject Classification:]
Primary:  05E10,  17B10;   Secondary:  81R05
\end{eremark}
\end{abstract}

\section{Introduction}   The general linear Lie algebra $\mathfrak{gl}(m)$
of $m \times m$ matrices over the complex numbers $\CC$ has a natural
action on the space $\VV$ of $m \times 1$ matrices  by matrix multiplication.
The symmetric group $\mathsf{S}_k$ acts by place permutations on  the
tensor power  $\VV^{\otimes k}$ and commutes with the action of $\mathfrak{gl}(m)$.  Schur  \cite{S1, S2}
exploited this commuting action
to accomplish the decomposition of $\VV^{\otimes k}$ into irreducible summands,
and by doing this for all $k$,  constructed all the irreducible polynomial representations of
$\mathfrak{gl}(m)$.  The commuting action determines an algebra epimorphism $\phi:  \CC \mathsf{S}_k \rightarrow \mathsf{End}_{\mathfrak{gl}(m)}(\VV^{\otimes k})$, which is an isomorphism
whenever $m\geq k$.    In the special case when $m=2$, the centralizer algebra
$\mathsf{End}_{\mathfrak{gl}(2)}(\VV^{\otimes k})$ is a quotient of the group
algebra $\CC \mathsf{S}_k$, which is the Temperley-Lieb algebra $\mathsf{TL}_k(2)$.

In their paper \cite{TL} on statistical mechanical models in
physics, H.N.V.  Temperley and E.H.  Lieb introduced what
are now known as the Temperley-Lieb algebras.  These algebras
play a prominent role in Jones's work (\cite {J1}, \cite{J3}) on
subfactors of von Neumann algebras. The Temperley-Lieb algebra  $\mathsf{TL}_k(\mathsf{q}+\mathsf{q}^{-1})$ also
can be viewed as a quotient of the Iwahori-Hecke algebra $\mathcal{H}_k(\qq)$
of type $\mathsf{A}$ (the above is the special case corresponding to the parameter $\mathsf{q}=1$),
and that realization led Jones \cite {J2} to
discover an invariant of knots and links (the Jones polynomial).
Kauffman's bracket model for the Jones polynomial has a tangle
theoretic interpretation using Temperley-Lieb algebras, and this
is closely connected with the recoupling theory of colored knots
and links and topological invariants of 3-manifolds \cite{KL}.

By comparison, the general linear Lie superalgebra $\mathfrak{gl}(m|n)$ over $\CC$
has a natural action on the space $\VV$ of $(m+n) \times 1$ matrices by matrix
multiplication.   The symmetric group $\mathsf{S}_k$ acts by {\it graded} place permutations on  the
tensor power  $\VV^{\otimes k}$  and commutes with the $\mathfrak{gl}(m|n)$-action.
As Berele and Regev show in \cite{BR},   the $\mathsf{S}_k$-action   yields an algebra epimorphism
$\phi:  \CC \mathsf{S}_k \rightarrow \mathsf{End}_{\mathfrak{gl}(m|n)}(\VV^{\otimes k})$, which is an isomorphism  whenever $(m+1)(n+1) > k$. This can be used to decompose $\VV^{\otimes k}$
into irreducible summands (see \cite{BL}),  and it enabled Berele and Regev to relate the
characters of the irreducible $\mathfrak{gl}(m|n)$-summands to  hook Schur functions.

In this work,  we examine the case $m=1=n$, which is analogous to the  $\mathfrak{gl}(2)$-case.
We show that the planar rook algebra $\CC \PS_{k-1}$ is the centralizer algebra
$ \mathsf{End}_{\mathfrak{gl}(1|1)}(\VV^{\otimes k})$ and use that insight to
give an explicit decomposition of  $\VV^{\otimes k}$ into irreducible summands.
The planar rook algebra, like its Temperley-Lieb counterpart, has a basis of diagrams.
It has a rich combinatorics and irreducible modules whose
dimensions are given by binomial coefficients (see \cite{FHH}).    The purpose of this
short note is to tie planar rook algebras to the
representation theory of the general linear Lie superalgebra $\gl$ and its quantum analogue.
In Section 6 of the paper, we consider the quantum enveloping algebra
$\UU_\qq(\gl)$ for $\qq$ not a root of unity  and tensor
powers of  its natural two-dimensional representation $\VV_\qq$.   Our main results
in that section are an explicit decomposition of $\VV_\qq^{\otimes k}$ into irreducible
modules for $\UU_\qq(\gl)$ and a proof that $ \mathsf{End}_{\UU_\qq(\gl)}(\VV_\qq^{\otimes k})$
is also isomorphic to $\CC \PS_{k-1}$.

In the quantum case, there is an algebra homomorphism from the
Iwahori-Hecke algebra $\mathcal H_k(\qq^2)$ of type $\mathsf A$  onto the centralizer algebra
$ \mathsf{End}_{\UU_\qq(\mathfrak{gl}(m|n))}(\VV_\qq^{\otimes k})$ for any $m$ and $n$,
where here $\VV_\qq$ is the natural $(m+n)$-dimensional representation
of $\UU_\qq(\mathfrak{gl}(m|n))$.  The irreducible $\UU_\qq(\mathfrak{gl}(m|n))$-summands
of $\VV_\qq^{\otimes k}$ have highest weights labeled by $(m|n)$-hook shape partitions of $k$,
and  the corresponding highest weight vectors were constructed in
\cite{M} using  Gyoja's $\qq$-analogue of Young symmetrizers in $\mathcal H_k(\qq^2)$.

In the special case when  $m=1,n=1$, Black, in recent work \cite{B}, applied Young's semi-normal representation
of $\mathcal H_k(\qq^2)$ on tableaux of hook shape (i.e. $(1|1)$-hook shape)  to
obtain an action of $\mathcal H_k(\qq^2)$ on sequences of length $k$ with components
$+$ or $-$.     The matrix units relative to  a basis indexed by
such sequences are described by
permutation diagrams with strands colored by $+$ or $-$.
Using the fact that there is a homomorphism $\mathsf{B}_k \rightarrow \mathcal H_k(\qq^2)$
from the braid group $\mathsf{B}_k$ on $k$ strands and identifying an element of
$\mathsf{B}_k$ with a link,  Black  calculates a state-sum formula for the Alexander polynomial
by taking traces of the action on sign sequences.

Our aim here is to determine the centralizer algebras $\mathsf{End}_{\mathfrak{gl}(1|1)}(\VV^{\otimes k})$ and
$ \mathsf{End}_{\UU_\qq(\mathfrak{gl}(1|1))}(\VV_\qq^{\otimes k})$
explicitly by realizing them as planar rook algebras  and to provide  an explicit decomposition of
tensor space in both cases;  problems not considered in \cite{B}.  In the final section of the paper we will explain how  the  colored permutation diagrams in \cite{B} are related to the planar rook diagrams in our approach.

\section{Planar Rook Algebras}

The rook monoid $\RR_n$ consists of all $n \times n$ matrices with at most one
1 in each row and column and 0 everywhere else.  Each matrix in $\RR_n$ corresponds
to an arrangement of nonattacking rooks on an $n \times n$ chessboard.   The permutation matrices inside $\RR_n$ give a copy of the  symmetric group $\mathsf{S}_n$.  Each matrix in
$\RR_n$ corresponds to an $n$-diagram;  that is,  a diagram with two rows with $n$
vertices in each and with edges that reflect the positions of the 1s:
$$\left(\begin{array}{ccccc}  0 & 0 & 0 & 0 & 0 \\
1 & 0 & 0 & 0 & 0 \\
0 & 0 & 0 & 0 & 1 \\
0 & 0 & 0 & 0 & 0 \\
0 & 1 & 0 & 0 & 0 \end{array}\right)
\ \  \longleftrightarrow \ \  {\beginpicture
\setcoordinatesystem units <0.4cm,0.2cm>
\setplotarea x from 1 to 7, y from -1.5 to 2.5
\put{$\bullet$} at  1 -1  \put{$\bullet$} at  1 2
\put{$\bullet$} at  2 -1  \put{$\bullet$} at  2 2
\put{$\bullet$} at  3 -1  \put{$\bullet$} at  3 2
\put{$\bullet$} at  4 -1  \put{$\bullet$} at  4 2
\put{$\bullet$} at  5 -1  \put{$\bullet$} at  5 2
\put{.} at 5.8 0
\plot 2 2  1  -1 /
\plot 5 2  2  -1 /
\plot 3 2  5  -1 /
\endpicture}$$
Matrix multiplication in  $\RR_n$ corresponds to concatenation of diagrams; so for example, if

$$d_1 =  {\beginpicture
\setcoordinatesystem units <0.4cm,0.2cm>
\setplotarea x from 1 to 7, y from -1.5 to 2.5
\put{$\bullet$} at  1 -1  \put{$\bullet$} at  1 2
\put{$\bullet$} at  2 -1  \put{$\bullet$} at  2 2
\put{$\bullet$} at  3 -1  \put{$\bullet$} at  3 2
\put{$\bullet$} at  4 -1  \put{$\bullet$} at  4 2
\put{$\bullet$} at  5 -1  \put{$\bullet$} at  5 2
\plot 2 2  1  -1 /
\plot 5 2  2  -1 /
\plot 3 2  5  -1 /
\endpicture} \quad \hbox{\rm and}  \quad  \ \  \
d_2 =  {\beginpicture
\setcoordinatesystem units <0.4cm,0.2cm>
\setplotarea x from 1 to 7, y from -1.5 to 2.5
\put{$\bullet$} at  1 -1  \put{$\bullet$} at  1 2
\put{$\bullet$} at  2 -1  \put{$\bullet$} at  2 2
\put{$\bullet$} at  3 -1  \put{$\bullet$} at  3 2
\put{$\bullet$} at  4 -1  \put{$\bullet$} at  4 2
\put{$\bullet$} at  5 -1  \put{$\bullet$} at  5 2
\plot 1 2  2  -1 /
\plot 2 2  4  -1 /
\put{,} at 5.8 0
\endpicture}$$
then

$$ {\beginpicture
\setcoordinatesystem units <0.4cm,0.2cm>
\setplotarea x from 1 to 7, y from -1.5 to 2.5
\put{$d_1d_2 \ \  =$} at -2 2
\put{$\bullet$} at  1 -1  \put{$\bullet$} at  1 2 \put{$\bullet$} at  1  5
\put{$\bullet$} at  2 -1  \put{$\bullet$} at  2 2  \put{$\bullet$} at  2 5
\put{$\bullet$} at  3 -1  \put{$\bullet$} at  3 2  \put{$\bullet$} at  3 5
\put{$\bullet$} at  4 -1  \put{$\bullet$} at  4 2  \put{$\bullet$} at  4 5
\put{$\bullet$} at  5 -1  \put{$\bullet$} at  5 2  \put{$\bullet$} at  5 5
\plot 1 2  2  -1 /
\plot 2 2  4  -1 /
\plot 2 5  1  2 /
\plot 5 5  2 2 /
\plot 3 5  5  2 /
\put{$ =$} at 7.5 2
\put{$\bullet$} at  10  .5  \put{$\bullet$} at  10 3.5
\put{$\bullet$} at  11  .5  \put{$\bullet$} at  11 3.5
\put{$\bullet$} at  12  .5  \put{$\bullet$} at  12 3.5
\put{$\bullet$} at  13 .5  \put{$\bullet$} at  13 3.5
\put{$\bullet$} at  14  .5  \put{$\bullet$} at  14 3.5
\put{.} at  14.8 1.5
\plot 11 3.5  11  .5 /
\plot 14 3.5  13  .5 /
\endpicture}$$

The rook algebra $\CC \RR_n$ has as a basis the elements of $\RR_n$,  and its unit element
is just  the diagram with $n$ vertical edges,  which corresponds to the identity matrix.  As shown in (\cite{Gr}, \cite{H}, \cite{HL}, \cite{So}),   the algebra $\CC \RR_n$ has many interesting combinatorial features and representation theoretic connections which come from the fact $\CC \RR_n$  generates the centralizer algebra  $\mathsf{End}_{\mathfrak {gl}(m)} ( (\CC^m \oplus \CC)^{\otimes n})$ for all $m \geq 1$.

The $n$-diagrams with no edge crossings form a submonoid $\PS_n$ of $\RR_n$, called the
{\it planar rook monoid}.  This can be regarded as the monoid of  all order-preserving partial
permutations of $\{1,2, \dots, n\}$,  see  \cite{R}.    The {\it planar rook algebra}  $\CC \PS_n$ has as a basis over $\CC$  the diagrams  in $\PS_n$ and multiplication that is  the linear extension
of the product in $\PS_n$.   By convention,  $\CC \PS_0 = \CC$.

To  count the number of diagrams in $\PS_n$ with exactly $\ell$ edges, we choose $\ell$ vertices from each row of the diagram, and then connect those vertices in the one and only one planar way (the first on the bottom  to the first on the top, the second to the second, etc.).  There
are exactly  ${n \choose \ell}^2$ such diagrams, and from this we see
$$\dim \CC \PS_n = | \PS_n |  = \sum_{\ell=0}^n    {n \choose \ell}^2 = \sum_{\ell=0}^n  {n \choose \ell} {n \choose n-\ell}  =  {{2n} \choose n}.$$
The last equality is a special case of the Chu-Vandermonde identity and can be
easily deduced from examining the coefficient of $z^n$  in the expansion of $(1 + z)^{n} (1 + z)^{n} = (1 + z)^{2n}$.

Given a diagram $d$ in $\PS_n$, let $\beta(d)$ (resp. $\tau(d)$)  be the vertices in the bottom (resp. top)  row with edges emanating from them.   Then $d$ can be regarded as an order-preserving bijective map $d: \beta(d) \rightarrow \tau(d)$  with domain $\beta(d)$
and image $\tau(d)$.    Thus, for the diagram
$$d =  {\beginpicture
\setcoordinatesystem units <0.4cm,0.2cm>
\setplotarea x from 1 to 7, y from -1.5 to 2.5
\put{$\bullet$} at  1 -1  \put{$\bullet$} at  1 2
\put{$\bullet$} at  2 -1  \put{$\bullet$} at  2 2
\put{$\bullet$} at  3 -1  \put{$\bullet$} at  3 2
\put{$\bullet$} at  4 -1  \put{$\bullet$} at  4 2
\put{$\bullet$} at  5 -1  \put{$\bullet$} at  5 2
\put{,} at  5.8 .2
\plot 2 2  1  -1 /
\plot 3 2  2  -1 /
\plot 4 2  5  -1 /
\endpicture}$$
$\beta(d) = \{1,2,5\}$ and $\tau(d) = \{2,3,4\}$.

When $\mathsf{s}$ and $\mathsf{t}$
are subsets of $\{1,\dots, n\}$ such that $| \mathsf{s} | = | \mathsf{t} |$,  we adopt
the notation $d_\mathsf{t}^\mathsf{s}$ for the diagram with $\beta(d_\mathsf{t}^\mathsf{s}) =
\mathsf{t}$ and $\tau(d_\mathsf{t}^\mathsf{s})  = \mathsf{s}$.

Following \cite{FHH}, we construct a $\CC \PS_n$-module $\M = \bigoplus_{\mathsf{s} \subseteq \{1,\dots, n\}} \CC m_{\mathsf{s}}$ with basis indexed by the subsets $\mathsf{s}$
of  $ \{1,\dots, n\}$,  where the $\CC \PS_n$-action is defined by

$$d m_{\mathsf{s}}  =  \begin{cases}    m_{d\mathsf{s}}   & \qquad \hbox{\rm if} \ \ \mathsf{s} \subseteq \beta(d) \\
0 & \qquad \hbox{\rm otherwise.} \end{cases} $$
Note that $d m_{\emptyset} = m_{\emptyset}$.
In particular, for the diagram $d$ pictured above,  $d m _{\{2,5\}} = m_{\{3,4\}}$,
while $d m _{\{4,5\}} = 0$.

\begin{prop} {\rm (\cite{FHH})}   Let $\M_\ell = \bigoplus_{|\mathsf{s}| = \ell}   \CC m_{\mathsf{s}} \subseteq \M$ for $\ell = 0,1,\dots, n$.   Then $\M_\ell$ is
an irreducible $\CC \PS_n$-submodule of $\M$ of dimension ${n \choose \ell}$
and  $\M = \bigoplus_{\ell=0}^n  \M_\ell$.  \end{prop}

\proof     Assume $\N$ is a nonzero submodule of $\M_\ell$ and
$0 \neq m' = \sum_{\mathsf{s}, |\mathsf{s}| = \ell}  \xi_\mathsf{s} m_{\mathsf s} \in \N$.
We suppose $\xi_{\mathsf{t}} \neq 0$.    Then
$d_\mathsf{t}^\mathsf{t} m' = \sum_{\mathsf{s}, |\mathsf{s}| = \ell}  \xi_\mathsf{s} d_\mathsf{t}^\mathsf{t} m_{\mathsf s} = \xi_{\mathsf{t}} m_{\mathsf{t}} \in \N$,  because in order for $d_\mathsf{t}^\mathsf{t}m_{\mathsf {s}}$ to be
nonzero,  $\mathsf{s} \subseteq \beta( d_\mathsf{t}^\mathsf{t}) = \mathsf{t}$ must hold;
but  $|\mathsf{s}| = |\mathsf{t}|$,  so  $\mathsf{s} = \mathsf{t}$ is forced.   Hence,
$m_{\mathsf{t}} \in \N$ and so is $d_\mathsf{t}^\mathsf{s} m_{\mathsf{t}} = m_{\mathsf{s}}
\in \N$ for all $\mathsf{s}$ with $|\mathsf{s}| = \ell$.  From this it follows that any nonzero
submodule $\N$ of $\M_\ell$ must equal $\M_\ell$, so $\M_\ell$ is irreducible.
The rest of the assertions follow readily.   \qed

\begin{prop}  {\rm (\cite[Thm.~3.2]{FHH})}   The modules $\M_\ell$ for $\ell = 0,1\dots, n$
form a complete set of irreducible $\CC \PS_n$-modules.
As a $\CC \PS_n$-module, $\CC \PS_n \cong \bigoplus_{\ell = 0}^n  {n \choose \ell} \M_\ell$.
Thus, the regular representation of $\CC \PS_n$ is completely reducible, and $\CC \PS_n$ is a semisimple associative algebra.    \end{prop}

\section{$\gl$}

The general linear  Lie superalgebra $\mathfrak{g} = \gl =\mathfrak g_0 \oplus \mathfrak g_1$ over $\CC$ consists
of all $2 \times 2$ complex matrices under the commutator product
$[a,b] = ab - (-1)^{\alpha\beta} ba$  for $a \in \mathfrak{g}_\alpha, b \in \mathfrak{g}_\beta$ and
$\alpha,\beta \in \mathbb Z/2 \mathbb Z$.    This product is bilinear and satisfies
$[b,a] = -(-1)^{\alpha\beta}[a,b]$  and the super-Jacobi identity
$[a,[b,c]] = [[a,b],c] + (-1)^{\alpha\beta}[b,[a,c]]$ for all $a \in \mathfrak{g}_\alpha, b \in \mathfrak{g}_\beta, c \in \mathfrak{g}$.

 The elements
$$e = \left ( \begin{array}{cc} 0 & 1 \\
0 & 0 \end{array} \right), \ \ f = \left ( \begin{array}{cc} 0 & 0 \\
1 & 0 \end{array} \right),  \ \ h_1 = \left ( \begin{array}{cc} 1 & 0 \\
0 & 0 \end{array} \right),  \ \ h_2 =
 \left ( \begin{array}{cc} 0 & 0 \\
0 & 1 \end{array} \right)$$
constitute a basis for $\gl$,  and  $\g_0 = \CC h_1 \oplus \CC h_2$  and $\g_1 = \CC e \oplus \CC f$.   It is easy to check that the
following relations hold:
$$[e,f] = ef+fe = h_1+h_2 = I, \ \
[h_1,e] = e, \ \  [h_2,e] = -e, \  \ [h_1,f] =- f, \ \ [h_2,f] = f.$$

In any Lie superalgebra,  an odd degree element $a$ satisfies  $a^2 = \frac{1}{2} [a,a]$ in the universal enveloping algebra.     In particular, in the universal enveloping algebra of $\gl$, we have
\begin{equation}\label{eq:squares} e^2 = \frac{1}{2} [e,e] = 0 = \frac{1}{2} [f,f] = f^2,\end{equation} because computing the commutators $[e,e]$ and $[f,f]$ in $\gl$
involves squaring  off-diagonal matrix units, which gives 0.

Each finite-dimensional  irreducible   $\gl$-module has a unique (up to
scalar multiple) highest weight vector; that is,  a nonzero vector $w$
such that  $e w = 0$ and
$h_1 w = m w$, $h_2 w = n w$ for some integers  $m,n$.  The module
is uniquely determined by the pair $[m,n]$,  and we denote it
by $\LL[m,n]$.

Here we consider the two-dimensional irreducible $\gl$-module $\VV = \VV_0 \oplus \VV_1$,
where $\VV_0 = \CC x$ and $\VV_1 = \CC y$,  and its $k$-fold tensor power
$\VV^{\otimes k}$  for all $k \geq 1$.  We can identify $\VV$ with the $2 \times 1$ matrices
over $\CC$, and $x$ and $y$ with the standard unit vectors

$$x = \left(\begin{array}{c}  1 \\ 0 \end{array}\right ) \qquad  y = \left( \begin{array}{c}  0 \\ 1\end{array}\right),$$
where the action of $\gl$ is given by matrix multiplication:
\begin{gather*}  ex = 0,   \quad  ey = x,  \qquad fx = y, \quad fy = 0, \\
h_1 x =  x,  \quad h_1y = 0 , \qquad  h_2 x = 0, \quad h_2 y = y. \end{gather*}
In the notation above,   $\VV = \LL[1,0]$ and
$x$ is the highest weight vector of weight $[1,0]$.

Now it is not difficult to argue that
$$\LL[m,n] \otimes \VV  =  \LL[m+1,n] \oplus \LL[m,n+1].$$
Starting with $\VV$ and iterating, we see that

\begin{align*} \VV^{\otimes 2} &=  \LL[2,0] \oplus \LL[1,1]\\
\VV^{\otimes 3} &=   \LL[3,0] \oplus 2 \LL[2,1] \oplus \LL[1,2] \\
\VV^{\otimes 4}  &=    \LL[4,0] \oplus 3 \LL[3,1] \oplus 3 \LL[2,2] \oplus \LL[1,3] \\
\vdots \ \  &  \\
\VV^{\otimes k} &= \bigoplus_{\ell=0}^{k-1} {{k-1} \choose \ell} \LL[\ell+1, k-1-\ell]  \end{align*}
The last line follows from an easy inductive argument which shows that  the multiplicity of $\LL[\ell+1,k-1-\ell]$ in $\VV^{\otimes k}$ is equal to the number of paths of length $k-1$ from  $\LL[1,0]  = \VV$  to $\LL[\ell+1,k-1-\ell]$  with $\ell$ steps in which 1 is added  to the first component  of the highest weight and $k-1-\ell$ steps with 1 added to the second component.

Since the $\gl$-module $\VV^{\otimes k}$ is completely reducible, by  the classical double-centralizer theory (see for example \cite[Secs. 3B and 68]{CR})  the centralizer algebra $\C_k = \mathsf{End}_{\gl} (\VV^{\otimes k})$  is a semisimple associative algebra over $\CC$, hence the sum of matrix blocks of size ${{k-1} \choose \ell} \times
{{k-1} \choose \ell}$,  and
$$\dim \C_k  = \sum_{\ell=0}^{k-1}  { {k-1} \choose \ell}^2   =   {{ 2(k-1)} \choose {k-1}}.$$
Thus,  $\C_k$ is isomorphic to the planar rook algebra  $\CC \PS_{k-1}$.   In
what follows, we want to exhibit the decomposition
of $\VV^{\otimes k}$ into irreducible $\gl$-summands and to demonstrate an explicit action of the
planar rook algebra $\CC\PS_{k-1}$ on $\VV^{\otimes k}$ which
commutes with the $\gl$-action.

\begin{section} {The Decomposition of $\VV^{\otimes k}$}  \end{section}

Our goal here is to display the decomposition of $\VV^{\otimes k}$ for
$\VV = \CC x \oplus \CC y$   into irreducible
$\gl$-submodules explicitly.   We use the fact that the action on a tensor product is given by
$$a (v \otimes w) = av \otimes w + (-1)^{\alpha \beta} v \otimes aw$$
for $a\in \gl$ of degree $\alpha$ and $v \in \VV$ of degree $\beta$.

Assume $\mathsf{s}$ is a subset of $\{1,\dots, k-1\}$ and define:
\begin{eqnarray}\label{eq:vp1}  u_{{\mathsf{s}}} &=&   u_1 \otimes \cdots \otimes u_{k-1} \\
v_{{\mathsf{s}}} &=&  e (y \otimes  u_{\mathsf{s}}) = x \otimes u_{\mathsf{s}} - y \otimes e u_{\mathsf{s}}, \label{eq:vp2}  \end{eqnarray}
where for $i = 1, \dots, k-1$,
$$u_i = \begin{cases}  x & \qquad \hbox{\rm if} \  i \in {\mathsf{s}} \\
y   & \qquad \hbox{\rm if} \  i \not \in {\mathsf{s}}. \end{cases}$$
For example, when $k = 4$ and ${\mathsf{s}} = \{1,3\}$, we have
$u_{{\mathsf{s}}} =  x \otimes y \otimes x$ and
$v_{{\mathsf{s}}}  = x \otimes x \otimes y \otimes x - y \otimes x \otimes x \otimes x.$
\medskip

\textbf{Claim 1.}    {\it For each subset   $\mathsf{s} \subseteq \{1,\dots, k-1\}$, the
vectors  $v_{\mathsf{s}}$  and  $fv_{\mathsf{s}}$ span
a two-dimensional irreducible $\gl$-submodule $\VV_{\mathsf{s}}$  of $\VV^{\otimes k}$  with highest weight vector $v_{\mathsf{s}}$ of weight $[|\mathsf{s}|+1, k-1-|\mathsf{s}|]$. \medskip

\proof Observe first that $v_{\mathsf{s}} \neq 0$ and
  $e v_{\mathsf{s}} = e^2 (y \otimes  u_{\mathsf{s}}) = 0$ by \eqref{eq:squares}.
Similarly,   $f(fv_{\mathsf{s}}) = f^2 v_{\mathsf{s}} =  0$.
Next, we have

$$e (fv_{\mathsf{s}}) = efv_{\mathsf{s}}  =   -fe v_{ \mathsf{s}} + I v_{\mathsf{s}}
= 0 + k v_{\mathsf{s}} \neq 0,$$
so in particular, $fv_{\mathsf{s}} \neq 0$.  Now
$$h_1 v_{\mathsf{s}} = h_1 e(y \otimes u_{\mathsf{s}})
= e h_1 (y \otimes u_{\mathsf{s}}) + e (y\otimes  u_{\mathsf{s}})
=  (|\mathsf{s}|+1) v_{\mathsf{s}}$$
and
$$h_2 v_{\mathsf{s}} = h_2 e (y \otimes  u_{ \mathsf{s}})
= e h_2 (y \otimes  u_{\mathsf{s}})   -e (y \otimes  u_{\mathsf{s}})
= (k-1-|\mathsf{s}|) v_{\mathsf{s}}.$$
Similarly,   $h_1 f v_{\mathsf{s}}= |\mathsf{s}| fv_{\mathsf{s}}$
and $h_2 f v_{\mathsf{s}}  = (k-|\mathsf{s}|) fv_{\mathsf{s}}$.

It follows that $\VV_{\mathsf{s}}$ is a two-dimensional submodule of $\VV^{\otimes k}$;
$v_{\mathsf{s}}$ is a highest weight vector of weight $[ |\mathsf{s}| +1, k-1- |\mathsf{s}|]$;
and $f v_{\mathsf{s}}$ is a lowest weight vector of weight $[\|\mathsf{s}|, k-|\mathsf{s}|]$.
It is easy to argue that $\VV_{\mathsf{s}}$
is irreducible,  since any nonzero submodule must decompose into weight spaces relative
to $h_1$ and $h_2$,
and so must contain either $v_{\mathsf{s}}$ or  $f v_{\mathsf{s}}$, hence both of them.
Therefore,  all the assertions in Claim 1 hold.
\medskip

\textbf {Claim 2.}  {\it  $\VV^{\otimes k} = \bigoplus_{\mathsf{s}} \VV_{\mathsf{s}}$, where $\mathsf{s}$ ranges over
the subsets of $\{1,\ldots, k-1\}$, affords a decomposition of $\VV^{\otimes k}$
into irreducible $\gl$-modules. }
\medskip

\proof If   $0 \neq  \VV_{\mathsf{t}} \cap \sum_{\mathsf{s} \neq
\mathsf{t}} \VV_{\mathsf{s}}$ for some subset $\mathsf{t}$,  then
$ \VV_{\mathsf{t}}$ must be contained in the right-hand sum by its
irreducibility.    Hence, $v_{\mathsf{t}}$ must be
a nonzero sum of highest weight vectors in the sum;  that is, it must be a linear
combination of the vectors  $ v_{\mathsf{s}}$ with $\mathsf{s} \neq
\mathsf{t}$  (and in fact,  with $|\mathsf{s}| = |\mathsf{t}|$).
However,  the terms $x \otimes u_{\mathsf{s}}$ in the vectors
$ v_{\mathsf{s}}$ in
\eqref{eq:vp2}
are all linearly independent.   This contradiction shows that the sum $\sum_{\mathsf{s}} \VV_{\mathsf{s}}$ is direct.
But $\dim \left(\bigoplus_{\mathsf{s}} \VV_{\mathsf{s}} \right) =  2^{k-1} \times 2 = 2^k = \dim (\VV^{\otimes k})$, so Claim 2 holds.

\begin{thm}  $ \displaystyle{\VV^{\otimes k} = \bigoplus_{\mathsf{s}} \VV_{\mathsf{s}} =
\bigoplus_{\ell = 0}^{k-1}   {{k-1} \choose \ell}   \LL[\ell+1, k-1-\ell]. } $ \end{thm}

\proof This is an immediate consequence of Claims 1 and  2.   \qed

\section{The action  of $\CC \PS_{k-1}$ on $\VV^{\otimes k}$}

It follows from the previous section that $\{ v_{\mathsf{s}}, f v_{\mathsf{s}} \mid \mathsf{s} \subseteq \{1,\dots, k-1\}\}$ is a basis for $\VV^{\otimes k}$.
Given a diagram $d \in \PS_{k-1}$,  we define an action of $d$ on $\VV^{\otimes k}$
by specifying its action on this basis according to

\begin{equation} \label{eq:dact} d v_{\mathsf{s}} =  v_{d \mathsf{s}},    \qquad  d (f v_{\mathsf{s}})  =  fv_{d \mathsf{s}} = f (d v_{\mathsf {s}}), \end{equation}
and then extend the action linearly  from $\PS_{k-1}$ to all of  $\CC \PS_{k-1}$.
Since $ d_1 (d_2 \mathsf{s}) = (d_1 d_2) \mathsf{s}$ holds for all $d_1, d_2 \in \PS_{k-1}$
and all subsets $\mathsf{s} \subseteq \{1,\dots, k-1\}$,  this makes $\VV^{\otimes k}$ into
a module for the algebra $\CC \PS_{k-1}$.     \bigskip

\textbf {Claim 3.}  {\it The $\CC \PS_{k-1}$-action and the $\gl$-action on $\VV^{\otimes k}$
commute.}

\medskip
\proof Observe that  $h_1d v_{\mathsf{s}} = h_1 v_{d\mathsf{s}} = |d \mathsf{s}| v_{d \mathsf{s}} = |\mathsf{s}| v_{d \mathsf{s}} = d h_1 v_{\mathsf{s}}$,    and  similarly for $h_2$.   As $de v_{\mathsf{s}} = 0 = e v_{d \mathsf{s}}
= e d v_{\mathsf{s}}$, it is apparent from this and \eqref{eq:dact} that $d$ commutes
with the action of $\gl$ on all the vectors $v_{\mathsf{s}}$.     It also commutes with the
action of $\gl$ on the vectors $fv_{\mathsf{s}}$.   For example,

$$de (fv_{\mathsf{s}}) = d(ef v_{\mathsf{s}})  = d(I v_{\mathsf{s}}) = k dv_{\mathsf{s}}= k v_{d\mathsf{s}}
= ef v_{d\mathsf{s}} = ef(d v_{\mathsf{s}}) = ed( fv_{\mathsf{s}}).$$
We leave the rest of the verifications to the reader.     \bigskip

 \textbf {Claim 4.}  {\it   If $p \in \CC \PS_{k-1}$ has the property that $p w = 0$
 for all $w \in \VV^{\otimes k}$, then $p = 0$.}

 \medskip
 \proof  We may suppose that $p = \sum_{\mathsf{s}, \mathsf{t}}
 \zeta_{\mathsf{t}}^{\mathsf {s}} d_{\mathsf{t}}^{\mathsf {s}}$
 where  $ \zeta_{\mathsf{t}}^{\mathsf {s}} \in \CC$ and  $\mathsf{s},\mathsf{t}$ are subsets of $\{1,\dots, k-1\}$
 with $|\mathsf{s} | = |\mathsf{t}|$.     We assume $\mathsf{t}'$ is chosen so
 $|\mathsf{t}'|$ is maximal among all subsets $\mathsf{t}$  with $\zeta_{\mathsf{t}}^{\mathsf{s}} \neq 0$ for some $\mathsf{s}$.     Then

 $$0 = p v_{\mathsf{t}'} =  \sum_{\mathsf{s}, |\mathsf{s}|  = |\mathsf{t}'|}
 \zeta_{\mathsf{t}'}^{\mathsf {s}}  d_{\mathsf{t}'}^{\mathsf {s}} v_{\mathsf{t}'} = \sum_{\mathsf{s}, |\mathsf{s}|  = |\mathsf{t}'|}   \zeta_{\mathsf{t}'}^{\mathsf {s}} v_{\mathsf{s}},$$
 \noindent  which by the independence of the vectors  $v_{\mathsf{s}}$ forces $ \zeta_{\mathsf{t}'}^{\mathsf {s}}
 = 0$, a contradiction.   Thus, $p = 0$.
 \bigskip

\begin{thm}  $\mathsf{End}_{\gl}(\VV^{\otimes k}) = \CC \PS_{k-1}$ \end{thm}

\proof Claims 3 and 4 imply that we may regard  $\CC \PS_{k-1}$ as a subalgebra of
$\mathsf{End}_{\gl}(\VV^{\otimes k})$.   However

$$\dim(\CC \PS_{k-1}) = {{2(k-1)} \choose {k-1}} = \dim (\mathsf{End}_{\gl}(\VV^{\otimes k})), $$
so that equality is forced.     \qed

\section{Tensor representations for quantum  $\gl$}

In this section we introduce the quantum enveloping algebra $\UU_\qq(\gl)$ and study
tensor powers of its natural two-dimensional module $\VV_\qq$.   When $\qq$ is not a root of unity
$\VV_\qq^{\otimes k}$ is completely reducible,  and we display its decomposition into irreducible
summands.   The centralizer algebra of the $\UU_\qq(\gl)$-action on $\VV_\qq^{\otimes k}$
is  shown to be the planar rook algebra $\CC \PS_{k-1}$ just as in the $\gl$ case.

Let $\PS$ be the free $\ZZ$-module given by $\PS = \ZZ \varepsilon_1 \oplus \ZZ \varepsilon_2$
and let $\PS^* = \ZZ h_1 \oplus \ZZ h_2$ be the dual module under the natural  bilinear
pairing $\langle \, , \,\rangle$  for which $\langle h_i, \varepsilon_j \rangle = \delta_{i,j}$.
We set $\alpha = \varepsilon_1 -\varepsilon_2$,  and let $\qq$ be a nonzero element of $\CC$.
We consider the unital Hopf algebra
$\UU_\qq = \uqgl$  over $\CC$  with generators
$E,F, \sigma, q^h (h \in \PS^*)$  and  relations

\begin{gather*} q^h = 1 \quad (\hbox{\rm for} \ \ h = 0),   \qquad q^{h} q^{h'} = q^{h+h'}, \\
q^hE = \qq^{\langle h,\alpha\rangle} Eq^h,  \qquad q^hF = \qq^{-\langle h,\alpha\rangle}Fq^h, \\    EF+FE = \frac{K - K^{-1}}{\qq - \qq^{-1}} \quad \hbox{\rm where} \ \ K = q^{h_1+h_2},  \\
E^2 = 0 = F^2, \\
\sigma E = -E \sigma, \quad  \sigma F = - F \sigma,  \quad \sigma q^h = q^h \sigma,   \quad \sigma^2 = 1,\end{gather*}
for all $h,h' \in \PS^*$.   The coproduct $\Delta$, counit $\mathsf{u}$,  and antipode $S$ on
$\UU_\qq$  are given by
\begin{gather*}   \Delta(E) =  E\otimes K^{-1} + \sigma \otimes E,  \quad \Delta(F) = F \otimes 1 + \sigma K \otimes F \\
\Delta(q^h) = q^h \otimes q^h,  \qquad  \Delta(\sigma) =
\sigma \otimes \sigma \\
\mathsf{u}(E) = 0  =  \mathsf{u}(F) = 0,  \quad \mathsf{u}(K) = 1 = \mathsf{u}(\sigma) \\
S(E) = -\sigma E K,  \quad S(F) = -\sigma K^{-1} F,  \quad S(q^h) = q^{-h},  \quad S(\sigma) = \sigma. \end{gather*}

\noindent The algebra $\UU_\qq$ is the quantum enveloping algebra  of the Lie superalgebra
$\gl$ studied in \cite{BKK}.

Here we consider the two-dimensional module $\VV_\qq = \CC x \oplus \CC y$  for $\UU_\qq$ and the
$k$-fold tensor power $\VV_\qq^{\otimes k}$ of $\VV_\qq$.       The action of $\UU_\qq$ on $\VV_\qq$
is given by
\begin{gather*}  Ex = 0,   \quad  Ey = x,  \qquad Fx = y, \quad Fy = 0, \\
q^h x = \qq^{\langle h,\varepsilon_1\rangle} x,  \quad q^h y = \qq^{\langle h,\varepsilon_2\rangle} y, \qquad  \sigma(x) = x, \quad \sigma(y) = -y. \end{gather*}
The coproduct gives the action of $\UU_\qq$ on a tensor product of any two $\UU_\qq$-modules.

Imitating the $\gl$ case,  we consider subsets $\mathsf{s}$ of $\{1,\dots, k-1\}$ and define

\begin{eqnarray}\label{eq:vp3} \quad u_{{\mathsf{s}}} &=&   u_1 \otimes \cdots \otimes u_{k-1} \\
\quad v_{{\mathsf{s}}} &=&  E (y \otimes  u_{\mathsf{s}}) = x \otimes K^{-1}u_{\mathsf{s}} - y \otimes Eu_{\mathsf{s}} = \qq^{1-k} x \otimes u_{\mathsf{s}} - y \otimes Eu_{\mathsf{s}}, \label{eq:vp4}  \end{eqnarray}
where for $i = 1, \dots, k-1$,
$$u_i = \begin{cases}  x & \qquad \hbox{\rm if} \  i \in {\mathsf{s}} \\
y   & \qquad \hbox{\rm if} \  i \not \in {\mathsf{s}}. \end{cases}$$
For example, when $k = 3$ and ${\mathsf{s}} = \{1\}$, we have
$u_{{\mathsf{s}}} =  x \otimes y$ and
$$v_{{\mathsf{s}}} = E(y \otimes x \otimes y) = \qq^{-2}x \otimes x \otimes y - y \otimes x \otimes x.$$

\textbf{Claim 5.}    {\it For each subset   $\mathsf{s} \subseteq \{1,\dots, k-1\}$, the
vectors  $v_{\mathsf{s}}$  and  $Fv_{\mathsf{s}}$ span
a two-dimensional irreducible $\UU_{\qq}$-submodule $\VV_{\qq,\mathsf{s}}$  of $\VV_\qq^{\otimes k}$  with highest weight vector $v_{\mathsf{s}}$ of highest weight
$[\qq^{|s|+1}, \qq^{k-1-|s|}]$ relative to $q^{h_1}, q^{h_2}$  whenever $\qq$ is not a root of unity.}  \medskip

\proof    First note that $Ev_{\mathsf{s}} = E^2(y \otimes u_{\mathsf{s}}) = 0$, while
$$q^{h_1} v_{\mathsf{s}} =  \qq^{|\mathsf{s}|+1} v_{\mathsf{s}},  \quad
q^{h_2} v_{\mathsf{s}} =  \qq^{k-1-|\mathsf{s}|}  v_{\mathsf{s}}, \quad  \sigma v_{\mathsf{s}} = (-1)^{k - |s|} v_{\mathsf{s}}.$$
Similarly  $F(Fv_{\mathsf{s}}) = F^2 v_{\mathsf{s}} = 0$,
$$q^{h_1}Fv_{\mathsf{s}} =  \qq^{|s|}v_{\mathsf{s}}, \quad q^{h_2} F v_{\mathsf{s}} = \qq^{k-|s|}v_{\mathsf{s}}, \quad \sigma F v_{\mathsf{s}} =
(-1)^{k+1-|s|} F v_{\mathsf{s}}.$$

Now

\begin{equation}\label{eq:ef} EF v_{\mathsf{s}} = -FEv_{\mathsf{s}} +  \left(\frac{K-K^{-1}}{\qq - \qq^{-1}}\right)v_{\mathsf{s}}
= [k] v_{\mathsf{s}},\end{equation}
where $[k]$ is the $\qq$-integer given by
$$[k] := \frac{\qq^k-\qq^{-k}}{\qq - \qq^{-1}}.$$
From this it is apparent that $\VV_{\qq, \mathsf{s}}$ is a two-dimensional $\UU_\qq$-module.  It is irreducible,  because any submodule must contain an eigenvector for $\sigma$
hence either $v_{\mathsf{s}}$ or $Fv_{\mathsf{s}}$.      But since
$[k] \neq 0$, it will then contain both vectors (see \eqref{eq:ef}).
\medskip

Let $\LL_\qq[m, n]$ denote the finite-dimensional irreducible $\UU_\qq$-module
with highest weight $[\qq^m, \qq^n]$ relative to $h_1, h_2$.   Then
$\VV_{\qq,\mathsf{s}} \cong  \LL_\qq[|s|+1, k-1-|s|]$, and as in the $\gl$ case
we have

\begin{thm}
  $ \displaystyle{\VV_\qq^{\otimes k} = \bigoplus_{\mathsf{s}} \VV_{\qq,\mathsf{s}} =
\bigoplus_{\ell = 0}^{k-1}   {{k-1} \choose \ell}   \LL_\qq[\ell+1, k-1-\ell], } $
 where $\mathsf{s}$ ranges over
the subsets of $\{1,\ldots, k-1\}$, affords a decomposition of $\VV_\qq^{\otimes k}$
into irreducible modules for $\UU_\qq = \UU_\qq(\gl)$ whenever $\qq$ is not a root of unity.
 \end{thm}

Just as in the $\gl$ case, we may define an action of the planar rook algebra $\CC \PS_{k-1}$on
$\VV_\qq^{k}$ by assigning for each diagram $d \in \PS_{k-1}$,

$$d v_{\mathsf{s}} =  v_{d\mathsf{s}} \qquad   dFv_{\mathsf{s}} = F  v_{d\mathsf{s}}  = Fdv_{\mathsf{s}}.$$

\begin{thm}  When $\qq$ is not a root of unity,  we have $\mathsf{End}_{\UU_q}(\VV_\qq^{\otimes k}) = \CC \PS_{k-1}$ for $\UU_\qq = \UU_\qq(\gl)$.  \end{thm}

\proof  The arguments to show that $\CC \PS_{k-1}$ commutes with the action of
$\UU_\qq$ on $\VV_\qq$ and to show that  $\CC \PS_{k-1}$ embeds into the
centralizer algebra $\mathsf{End}_{\UU_q}(\VV_\qq^{\otimes k})$ are virtually identical
to the $\gl$ case.   Here is one sample computation to show that a diagram $d \in \PS_{k-1}$
commutes with $E$ on a vector $F  v_{\mathsf{s}}$:

\begin{align*} (d E)F  v_{\mathsf{s}} &= d (EF) v_{\mathsf{s}}  = -d(FE) v_{\mathsf{s}}
+ d \left(\frac{K-K^{-1}}{\qq - \qq^{-1}}\right)v_{\mathsf{s}} \\
&= [k] d v_{\mathsf{s}}  =  [k] v_{d\mathsf{s}}  =  \left(\frac{K-K^{-1}}{\qq - \qq^{-1}}\right)v_{d\mathsf{s}}  \\
&= EF v_{d\mathsf{s}}   = EF d v_{\mathsf{s}}   = EdF v_{\mathsf{s}}  \\
&= (Ed) Fv_{\mathsf{s}}.  \end{align*}
The remaining verifications are left to the reader.   \qed

\section{Connections with representations of $\mathcal H_k(\qq^2)$ }

It is known (see \cite{M}, for example) that there is an algebra epimorphism
$\varphi: \mathcal H_k(\qq^2) \rightarrow \mathsf{End}_{\UU_\qq(\gl))}(\VV_\qq^{\otimes k})$
giving a representation of the Iwahori-Hecke algebra $ \mathcal H_k(\qq^2)$ on $\VV_\qq^{\otimes k}$. The quotient $\mathcal H_k(\qq^2)/\ker \varphi$ is semisimple and decomposes into matrix
blocks indexed by the partitions of $k$ of hook shape.   Black \cite{B} described an action of this quotient algebra on the span of vectors  $x_\eta$
indexed by sequences $\eta = (\eta_1, \dots, \eta_k)$ with $\eta_j = \pm$ for $j=1,\dots,k$ such that $\eta_1 = +$.    The sequences $\eta$ having  $\ell$ coordinates equal to $-$ are in one-to-one correspondence with the
standard tableaux whose shape is given by the partition $(k-\ell, 1^{\ell})$, and the action of
$\mathcal H_k(\qq^2)$ on the vectors $x_\eta$ indexed by those sequences is  derived from Young's semi-normal representation
on the corresponding standard tableaux.

Now the matrix units relative to the basis $x_{\eta}$ are in bijection with
the permutation diagrams with strands colored  $+$ and $-$  such that
no two strands of the same color cross, and the first strand goes directly down and has color $+$.
For example,  the matrix unit $\mathsf{E}_{\vartheta, \eta}$  labeled by the sequences $\vartheta = (+,+,+,-,-,-)$ and  $\eta = (+,-,+,-,+,-)$  corresponds to the permutation diagram}

$$
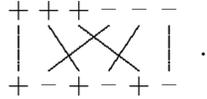
{\beginpicture
\setcoordinatesystem units <0.4cm,0.2cm>
\setplotarea x from 0 to 6, y from -1.5 to 2.5
\plot 0 2  0  -1 /
\plot 1 2  2  -1 /
\plot 2 2  4  -1 /
\plot 3 2  1  -1 /
\plot 4 2  3  -1 /
\plot 5 2  5  -1 /
\put{+} at 0  -2
\put{+} at 2  -2
\put{+} at 4  -2
\put{--} at 1  -2
\put{--} at 3  -2
\put{--} at 5  -2
\put{+} at 0  3
\put{+} at 2  3
\put{--} at 4  3
\put{+} at 1  3
\put{--} at 3  3
\put{--} at 5  3

\endpicture}.
$$

There is a homomorphism $\psi:  \mathsf{B}_k \rightarrow \mathcal H_k(\qq^2)$ from
the braid group $\mathsf{B}_k$ on  $k$-strands,  and so the
action above determines an action of $\mathsf{B}_k$
on the space spanned by the vectors   $x_\eta$. Taking  traces gives the Alexander
polynomial of the link associated to a braid group element.

The correspondence between the colored permutation diagrams and the elements
of the planar rook monoid $\mathsf{P}_{k-1}$ can be obtained by ignoring the first
strand and deleting the strands colored $-$.  The diagram above
corresponds to  the element $\mathsf{X}_d = \sum_{d' \subseteq d}  (-1)^{|d \setminus d'|} d' \in
\mathbb C \mathsf{P}_{k-1}$
in \cite{FHH}  for

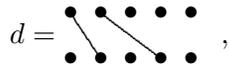
$$
d =  {\beginpicture
\setcoordinatesystem units <0.4cm,0.2cm>
\setplotarea x from 1 to 6, y from -1.5 to 2.5
\put{$\bullet$} at  1 -1  \put{$\bullet$} at  1 2
\put{$\bullet$} at  2 -1  \put{$\bullet$} at  2 2
\put{$\bullet$} at  3 -1  \put{$\bullet$} at  3 2
\put{$\bullet$} at  4 -1  \put{$\bullet$} at  4 2
\put{$\bullet$} at  5 -1  \put{$\bullet$} at  5 2
\plot 1 2  2  -1 /
\plot 2 2  4  -1 /
\endpicture},
$$
where $d'$ ranges over all the diagrams obtained from $d$ by deleting
edges, and $| d \setminus d'|$ is the number of edges in $d$ minus the number of edges
in $d'$.    If $d_1$ and $d_2$ are two diagrams in $\mathsf{P}_{k-1}$,  then by \cite[Prop.~3.3]{FHH} we have

$$\mathsf{X}_{d_1} \mathsf{X}_{d_2}   = \delta_{\beta(d_1), \tau(d_2)} \mathsf{X}_{d_3},$$
where $d_3$ is the diagram with top row $\tau(d_1)$ and bottom row $\beta(d_2)$.   Therefore, the correspondence
between the matrix units $\mathsf{E}_{\vartheta, \eta}$ and the matrix units $\mathsf{X}_d$ determines
an algebra isomorphism between the algebra having basis the colored permutation diagrams on sequences of length $k$  and the planar rook algebra $\mathbb C  \mathsf{P}_{k-1}.$


\begin{thebibliography}{FHH}

 \bibitem[BKK]{BKK}  G.~Benkart, S.-J. Kang, and M.~Kashiwara, Crystal bases for the
quantum superalgebra $U_q(gl(m,n))$, {\it J.~Amer.  Math.  Soc.}  {\bf 13} (2000), 295-313.



\bibitem[BL]{BL} G.~Benkart and C.~Lee (Shader),  Stability in modules for general linear Lie superalgebras,  {\it Nova J.  of Algebra and Geometry} {\bf 2} (1994),
383-409.


\bibitem[BR]{BR}  A. Berele and A. Regev, Hook Young diagrams with applications to combinatorics and to representations of Lie superalgebras,  {\it Adv. in Math.} \textbf{64}
(1987), 118-175.

\bibitem[B]{B} S. Black,  A state-sum formula for the Alexander polynomial, arXiv: 1002.4860v2.

\bibitem[CR]{CR} C.~Curtis and I.~Reiner, {\em Methods of Representation Theory -- With Applications to Finite Groups and Orders}, Pure
and Applied Mathematics, vols. I and II, Wiley \& Sons, Inc., New York, 1987.

\bibitem[FHH]{FHH}D.~Flath, T.~Halverson, and K.~Herbig, The planar rook
algebra and Pascal's triangle, {\it Enseign.~Math.}~(2) {\bf 54} (2008), 1-16.

\bibitem[Gr]{Gr}C.~Grood,  A Specht module analog for the rook monoid, {\it Electronic
J. Combin.} \textbf{R2 9} (2002).

\bibitem[H]{H}T.~Halverson,  Representations of the $q$-rook monoid, {\it J. Algebra}
 \textbf{273} (2004), 227-251.

 \bibitem[HL]{HL}T.~Halverson and T.~Lewandowski,  RSK insertion for  set partitions and diagram algebras, {\it Electronic
J. Combin.} \textbf{11 R24} (2004/06).

\bibitem[J1]{J1}
V.F.R.  Jones, Index for subfactors,  \emph{Invent. Math.}
 \textbf{72} (1983), 1--25.

 \bibitem[J2]{J2}
V.F.R.  Jones,  A polynomial invariant for knots via von Neumann
algebras, \emph{Bull.  Amer.  Math.  Soc.}  \textbf{12} (1985),
103--111.

\bibitem[J3]{J3}V.F.R.  Jones,
\emph{Subfactors and Knots}, CBMS Regional Conference Series in Mathematics, \textbf{80}. Published for the Conference Board of the Mathematical Sciences, Washington, DC, by the American Mathematical Society, Providence, RI, 1991.

\bibitem[KL]{KL}
L.H.~Kauffman and S.L.~Lins, Temperley-Lieb Recoupling
Theory and Invariants of 3-Manifolds, \emph{Annals of Mathematics
Studies} \textbf{134} Princeton Univ.  Press, Princeton 1994.

\bibitem[M]{M} D.~Moon, Highest weight vectors of irreducible representations of the
quantum superalgebra $\mathfrak U_\qq(gl(m,n))$,  {\it J.~Korean Math. Soc.}
\textbf{40} (2003), no. 1, 1--28.

\bibitem[R]{R}L.~Renner,  {\it Linear Algebraic Monoids}, Encyclopedia of Mathematical Sciences \textbf{134} Springer-Verlag, Berlin,  2005.

\bibitem[S1]{S1} I. Schur, {\it \"Uber eine Klasse von Matrixen die sich einer gegebenen
Matrix zuordnen lassen}, Thesis Berlin (1901), reprinted in I. Schur,
Gesammelte Abhandlungen I Springer, Berlin, (1973), 1-70.

\bibitem[S2]{S2} I. Schur, {\it \"Uber die rationalen Darstellungen der allgemeinen
linearen Gruppe} (1927), reprinted in I. Schur, Gesammelte Abhandlungen III
Springer, Berlin, (1973), 68-85.

\bibitem[So]{So}L.~Solomon, Representations of the rook monoid, {\it J. Algebra}
 \textbf{256} (2002), 309-342.

\bibitem[TL]{TL}
H.N.V.~Temperley and E.H.~Lieb, Relations between the
``percolation'' and  ``colouring'' problem and other
graph-theoretical problems associated with regular planar
lattices: some exact results for the ``percolation'' problem,
\emph{Proc.  Roy.  Soc.  London Ser.  A} \textbf{322}  (1971), 251--280.

\end{thebibliography}
\end{document}